\documentclass{article}

\usepackage{amsmath,amssymb,graphicx,subcaption,tikz,wrapfig}
\usepackage[top=1in,bottom=1in,left=1in,right=1in]{geometry}
\usepackage[numbers]{natbib}

\newcommand{\pd}[2]{\frac{\partial #1}{\partial #2}}
\newcommand{\dd}[2]{\frac{\mathrm d#1}{\mathrm d#2}}
\newcommand{\R}{\mathbb R}
\newcommand{\abs}[1]{\left \lvert #1 \right \rvert}
\newcommand{\dphi}[2]{\phi_{#1}^{#2}}
\newcommand{\ext}[1]{\underset{#1}{\text{ext}}}

\graphicspath{{/}}

\title{Modeling Environmental Crime in Protected Areas Using the Level Set Method}

\author{D. J. Arnold\thanks{Department of Mathematics, University of California Los Angeles.}
\and D. Fernandez\footnotemark[1]
\and R. Jia\footnotemark[1]
\and C. Parkinson\footnotemark[1]
\and D. Tonne\thanks{Department of Mathematics, California State University Long Beach.}
\and Y. Yaniv\thanks{Department of Mathematics, University of Maryland.}
\and A. L. Bertozzi\footnotemark[1]
\and S. J. Osher\footnotemark[1]}

\begin{document}

\maketitle

\begin{abstract}
National parks often serve as hotspots for environmental crime such as illegal deforestation and animal poaching. Previous attempts to model environmental crime were either discrete and network-based or required very restrictive assumptions on the geometry of the protected region and made heavy use of radial symmetry. We formulate a level set method to track criminals inside a protected region which uses real elevation data to determine speed of travel, does not require any assumptions of symmetry, and can be applied to regions of arbitrary shape. In doing so, we design a Hamilton-Jacobi equation to describe movement of criminals while also incorporating the effects of patrollers who attempt to deter the crime. We discuss the numerical schemes that we use to solve this Hamilton-Jacobi equation. Finally, we apply our method to Yosemite National Park and Kangaroo Island, Australia and design practical patrol strategies with the goal of minimizing the area that is affected by criminal activity.
\end{abstract}


\section{Introduction}
Environmental crime in protected national parks is a concern to authorities around the world. National parks often serve as hotspots for illegal logging and animal poaching, henceforth referred to as illegal extraction. In recent years, scientists have aided law enforcement agencies in the prevention of such crimes. One way in which they do this is by building models to describe deforestation, track animal movement, and predict adversarial behavior of criminals. \citet{leader1993policies} examined a case study of poaching in the Luangwa Valley, Zambia and concluded that increasing detection rates was a more effective deterrent to environmental crime than increasing severity of punishment. This means that models that can help improve detection rates of environmental criminals are a useful tool for patrollers of protected areas.

Authorities employ different strategies to patrol protected areas and detect illegal extractors, and research into the effectiveness of these strategies is of major interest. Researchers are challenged to find efficient and practical patrol strategies to police vast areas with limited resources. 
An important model of illegal extraction was given by \citet{Albers}. They describe illegal extraction as a continuous spatial game between patrol units and extractors, constructing a model where extractors maximize their expected profit by trading off costs of penetrating further into a protected region (increasing effort required and risk of capture) against increased benefits of extracting further from the boundary of the protected region. Albers' model was formulated as a Stackelberg game, an adversarial game with a defender (patrollers) and many attackers (extractors) with perfect information about the defender's strategy. Albers gave some qualitative results in a very simplified situation where the protected region was circular, and all quantities were radially symmetric. \citet{johnson2012patrol} proved an additional result regarding optimal patrol strategies for this simple case.

Subsequently, models for illegal extraction in protected regions using discrete rather than continuous methods have been developed by Fang et al. \citep{fang2013optimal,fang2017paws} and Kar et al. \citep{kar2015game,Intercept2017}. These models incorporate different methods for modeling human behavior and different ways to treat the protected region. \Citet{kar2015game} used repeated Stackelberg games to attempt to understand the evolution of attackers strategies, but did not consider realistic spatial domains. \Citet{fang2017paws} developed a model including detailed terrain and spatial information by describing the protected region as a series of nodes connected by edges corresponding to natural pathways through the region (for example along rivers or walking trails). The PAWS algorithm developed by \citet{fang2017paws} has been deployed in Queen Elizabeth National Park (QENP) in Uganda, and in Malaysian forests. Taking another approach, \citet{Intercept2017} used machine learning techniques to construct a model for predicting poaching attacks, again deploying the model for a field test in QENP, Uganda.

Our goal is to generalize the model of \citet{Albers} so that it is applicable to realistic protected regions, can include the effects of terrain and geometry, and can be executed for real protected regions. The most important mathematical technique that enables this extension is the level set method of \citet{osherSethian1988}, which provides a simple way to track moving fronts through the protected region. Indeed, the level set method and Hamilton-Jacobi equations are now used extensively in modeling travel in a variety of contexts. Sethian and Vladimirsky \cite{SethVlad2,SethVlad1} discuss a level set method for optimal travel on manifolds. In application, \citet{Dubins}, and later \citet{Takei}, modeled movement of simple cars through terrain with obstacles using methods from control theory which include Hamilton-Jacobi-Bellman equations. Tomlin has used similar methods to determine reachable sets for aircraft autolanders \cite{TomlinAircraft}, and to model human movement in adversarial reach-avoid games \cite{TomlinReachAvoid}. Thus there is a precedent for using Hamilton-Jacobi type equations to model movement throughout domains. However, as discussed above, attempts to model environmental crime have either required overly restrictive assumptions regarding radially symmetry and geometry of the protected region, or have been discrete in nature. This work serves to fill a gap in the literature. We suggest a continuum model for environmental crime which is still able to account for realistic terrain information, and bridges the divide between environmental crime modeling and Hamilton-Jacobi formulations for modeling human movement. 

This paper is structured as follows. In section \ref{sec:albersmodel}, we describe the model of \citet{Albers}. In section \ref{sec:levelset}, we introduce the level set method, describe the numerical methods used to solve it, and show how it can be applied to model movement through regions with terrain. In section \ref{sec:model}, we describe our model for illegal extraction from protected regions, and provide an algorithm to calculate the regions through which extractors will pass. In section \ref{sec:results}, we present results and discuss their implications towards both patrol strategies and our modeling approach. Finally, in section \ref{sec:conclusion}, we make conclusions and identify potential directions for further research. 



\section{The Illegal Logging Model of \citet{Albers}}\label{sec:albersmodel}
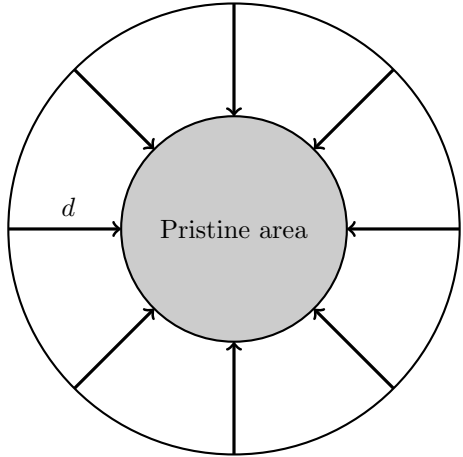
\begin{figure}
	\centering
\begin{tikzpicture}
\draw [thick,fill=black!20] (0,0) circle (1.5);
\draw [thick] (0,0) circle (3);
\draw [very thick,<-] (0,1.5)--(0,3);
\draw [very thick,<-] (0,-1.5)--(0,-3);
\draw [very thick,<-] (1.5,0)--(3,0);
\draw [very thick,<-] (-1.5,0)--(-3,0);
\draw [very thick,<-] (1.0607,1.0607) -- (2.1213,2.1213);
\draw [very thick,<-] (-1.0607,-1.0607) -- (-2.1213,-2.1213);
\draw [very thick,<-] (-1.0607,1.0607) -- (-2.1213,2.1213);
\draw [very thick,<-] (1.0607,-1.0607) -- (2.1213,-2.1213);
\node at (-2.2,0.3) {$d$};
\node at (0,0) {Pristine area};
\end{tikzpicture}
    \caption{Illustration of the \citet{Albers} model, where extraction occurs at depth $d$ in a circular region. Reproduction of Figure 2 from \citet{johnson2012patrol}.}
    \label{fig:albersmodel}
\end{figure}

\Citet{Albers} presented a model for illegal extraction from protected regions that includes spatial effects. They use a Stackelberg game in a circular region between defenders and extractors. In their model, an extractor gains benefit and incurs cost based on the depth $d$ into the region where they choose to extract. The benefit is associated with the extraction point, and is modeled as a concave increasing function $B(d)$. Cost is based on the trip distance, and is modeled as a quasiconvex increasing function $C(d)$, but no specific functions are given. Then the extractor's profit associated with extracting at depth $d$ is \begin{equation} \label{eq:AlbersProfitNoPatrol} P(d) = B(d) - C(d). \end{equation} If the region is patrolled with density $\psi(d)$, the probability that the extractor is detected is the cumulative probability $\Psi(d)=\int^d_0\psi(r)dr$ over the trip out of the protected region. If the illegal extractor is caught, any resources they've extracted are confiscated and they must leave empty-handed. The extractor's expected profit can therefore be written as \begin{equation} \label{eq:AlbersProfit} P(d) = (1 - \Psi(d))B(d) - C(d). \end{equation} Given the extractors have perfect information of the patrol strategy, they know the function $\Psi(d)$, and hence their task is to find the optimal distance $d^*$ to penetrate into the protected region. Simultaneously, the patrollers task is to pick the patrol strategy $\psi(d)$ that minimizes $d^*$. Figure \ref{fig:albersmodel} illustrates the model, showing the pristine region where extractors are never present and the outer region through which they pass while travelling into or out of the protected region.

\section{The Level-Set Method}\label{sec:levelset}
The level-set method is a powerful technique that we use to model human movement through complex domains. Originally developed by \citet{osherSethian1988}, the level set method models fronts which propagate with a prescribed velocity. The front is embedded as the zero level set of an auxiliary function, which is then evolved according to a partial differential equation. We describe a simple version of the method in $\R^2$, sufficient for our needs. Suppose we are given a simple, closed curve (or a collection of non-intersecting, simple, closed curves) denoted by $\Gamma$. To evolve $\Gamma$ via level set flow, we begin by finding a Lipschitz continuous function $\phi_0 : \R^2 \to \R$ such that $\phi_0$ is positive inside $\Gamma$ and negative outside $\Gamma$, and thus $\Gamma = \{ x \in \R^2  \, : \, \phi_0(x) = 0\}$. Next, we evolve $\phi:\R^2\times[0,\infty)\to\R$ using the PDE \begin{equation} \label{eq:genLevelSet} \begin{aligned} \phi_t + v(x) \abs{\nabla \phi} = 0,& \\ \phi(x,0) = \phi_0(x),& \end{aligned} \end{equation} for some non-negative velocity function $v$. Define $\Gamma(t) = \{x \in \mathbb R^2 \, : \, \phi(x,t) = 0\}$. This represents an evolution of $\Gamma = \Gamma(0)$. Indeed, we can re-write equation \ref{eq:genLevelSet} as $\phi_t + v(x) \left( \frac{\nabla \phi}{\abs{\nabla \phi}} \right) \cdot \nabla \phi = 0.$ Note that $\nabla \phi/ \abs{\nabla \phi}$ is a unit vector which is normal to the level contour $\Gamma(t)$. Thus, locally, this equation models advection in the direction normal to $\Gamma(t)$ with velocity $v(x)$. This causes $\Gamma(t)$ to deform with normal velocity $v(x)$. With this in mind, one can view $\Gamma(t)$ as the set of points which could be reached after walking inward from the boundary for a time $t$ with normal velocity $v(x)$. For example, if $v(x) \equiv 1$, then for $t > 0$, the curve $\Gamma(t)$ is the set of points a distance $t$ from the original curve.

\begin{figure*}\centering
\begin{subfigure}[t]{0.35\textwidth}
\includegraphics[width=\textwidth]{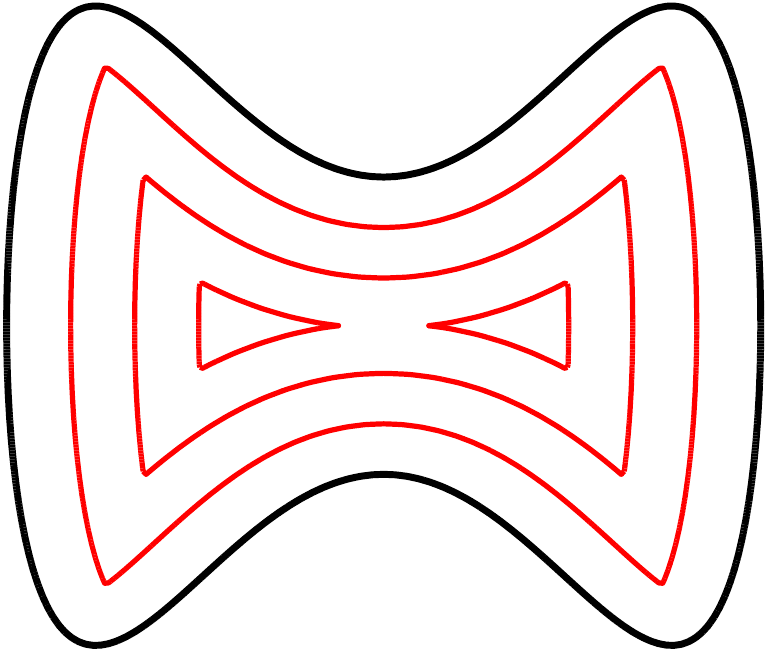}
\caption{Contours of equal travel time (red) from the boundary of a region (black).}
\label{fig:distEx}
\end{subfigure}~
\begin{subfigure}[t]{0.45\textwidth}
\includegraphics[width=\textwidth]{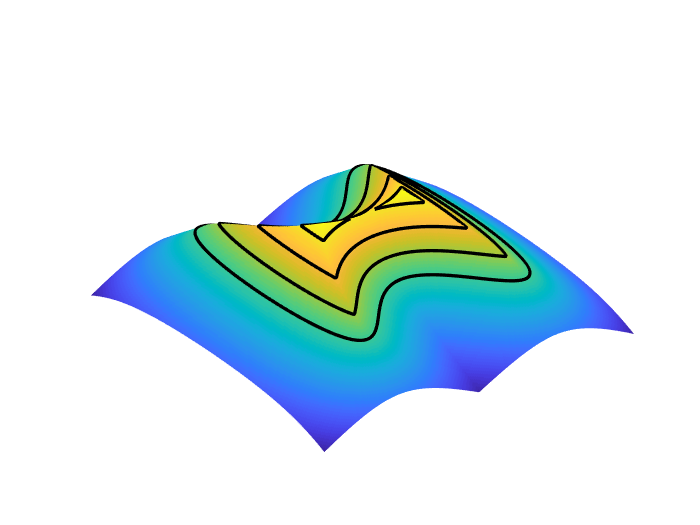}
\caption{The auxiliary function $\phi$ with level sets (black).}
\label{fig:initFun}
\end{subfigure}
\caption{Using the level set method to find contours of equal travel time.}
\label{fig:levelSetExample}
\end{figure*}

This procedure is shown in figure \ref{fig:levelSetExample}. Figure \ref{fig:distEx} shows the boundary of a region defined by $(x(\theta),y(\theta)) = (\cos(\theta) ,\sin(\theta)+\sin(3\theta)/2)$ and a series of contours of points accessible after traveling a certain time from the boundary. Figure \ref{fig:initFun} shows the graph of the level set function and with the same level sets overlaid on the surface. This example also displays one of the strengths of the level set method. Depending on the original curve, the the level contours may break apart into disconnected pieces or merge into into a single piece. More primitive methods for curve evolution (many of which involve tracking points along the curve) have trouble handling these changes in topology but the level set method accounts for them with no special considerations. For more information regarding level set methods, see \citet{osher2003level}.

\section{Illegal Deforestation}\label{sec:model}
Our model for the actions of illegal extractors committing crimes in protected regions is based on \citet{Albers}, as described in section \ref{sec:albersmodel}. They derived a game-theoretic model for illegal deforestation in an idealized setting, assuming the protected region was circular to significantly simplify the problem, and neglecting the effects of terrain. We wish to relax these restrictions which leads to several significant difficulties.

%

We are interested in real national parks and hence the protected region will not be symmetric and the benefit and cost functions could be arbitrary. In \citet{Albers}, the profit is maximized on a ring with some radius and extractors will enter from any point on the boundary of the region to travel to the maximum profit ring. Without radial symmetry however, the profit will, in general, be maximized at an isolated point rather than along a closed curve, and correspondingly there will be a single optimal route from this point to the boundary of the protected region. Hence illegal extractors only ever occupy a set of points of measure zero. In our model, we resolve this issue by assuming that extractors will tolerate extracting anywhere where the profit is close to the maximum possible profit.

The cost function of \citet{Albers} does not generalize easily to irregular geometries. Instead of using a predefined formula based on distance, we constructively determine the cost by considering the effects that would detract from the final profit with the help of the level set method. Given terrain and elevation information for a real protected park, the paths taken by extractors will not be simple straight lines. Therefore it is also necessary to find the paths extractors will follow to exit of the protected region. In the next section we describe an algorithm to find the pristine and non-pristine regions in a protected area with arbitrary shape, terrain, benefit distribution and patrol strategy.

\subsection{Problem Description}
Let $\Omega$  represent the area that needs to be protected and let $\psi(x,y)$ be the patrol density function measuring the likelihood of an extractor being caught at position $(x,y)$. The patrol budget $E$ is defined as $E=\int_{\Omega} \psi(x,y) \,\mathrm dx\,\mathrm dy$. The budget measures how many patrolling resources can be allocated to protect the domain $\Omega$. Effectively, the budget scales $\psi$.
 Extracting at a position $(x,y)$ gives the extractor an amount of benefit $B(x,y)$. For example, the benefit could depend on the quantity, value, or species of trees. We make no assumptions about $B$ except that it is known. $C(x,y)$ is the expected cost associated with extracting at position $(x,y)$. The cost is based on two factors, the effort involved in traveling from the extraction point to the boundary of the protected region, and the risk of being caught by patrols whilst inside the domain. The further into the protected the region the extractors penetrate, the higher the cost as they must expend more time/energy traveling, and are more likely to be captured by patrollers. 

This highlights a difference between our model and \citet{Albers}: we include the effect of the patrol directly in the cost function, whereas they include it as a modification to the benefit. In our model, the profit an extractor expects from extracting at a point $(x,y)$ in the protected region is simply $P(x,y)=B(x,y)-C(x,y)$. We assume the extractors accept a position to extract if the profit is within some tolerance of the maximum obtainable profit. Finally, we determine the paths that extractors will take when leaving the protected area (where they can be captured by the patrol). We assume the extractors will always take an optimal path from the extraction point to the boundary, and can only be caught after they have finished extracting and started heading back to the boundary (since no crime has been committed before they have extracted).

\subsection{Calculating the Profit Function}

Since the benefit function is known in advance, it remains to calculate the expected cost function. As explained in the previous section, the cost function depends on many factors, including the optimal route to a given point. The level set method avoids most of the difficulties one might expect when calculating optimal paths from all points inside the protected region to the boundary. We use the level set method to find contours of equal cost in the protected region. The level set method implicitly uses the optimal path between the boundary and an interior point, and can be efficiently implemented. 

The expected cost associated with extracting at a given point $ (x_0,y_0) $ is calculated with the following level set equation, using the boundary of the domain $\Omega$ as the zero-level set of the function $\phi(x,y,0)$:
\begin{equation} \label{eq:costLevelSet}
\pd{\phi}{t} = -\frac{1}{1/v(x,y)+\alpha\psi(x,y)B(x_0,y_0)} \abs{\nabla \phi}.
\end{equation}

\begin{wrapfigure}{r}{0.37\textwidth}
\centering
\includegraphics[width=.37\textwidth, trim = 25 10 50 20, clip]{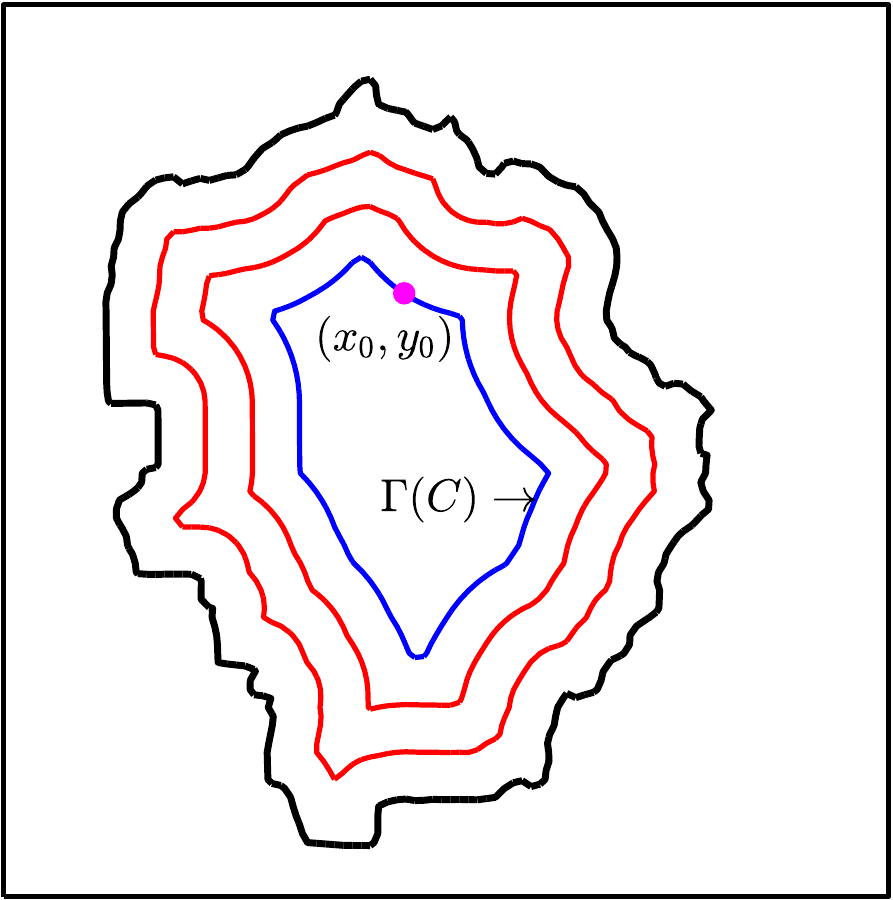}
\caption{Evolve level sets in from the boundary until the $C > 0$ such that $(x_0,y_0) \in \Gamma(C)$.} 
\label{fig:demonstrateCost}
\end{wrapfigure}

\noindent The cost $C(x_0,y_0)$ is implicitly defined by $\phi(x_0,y_0,C(x_0,y_0)) = 0$. Intuitively, the point $(x_0,y_0)$ lies somewhere inside the curve $\Gamma = \partial \Omega$. When we evolve the boundary using \ref{eq:costLevelSet}, the zero level-sets $\Gamma(t) = \{\phi(x,y,t) = 0 \}$ represent the set of points which require equal cost to reach from the boundary. For some $C > 0$, we will have $(x_0,y_0) \in \Gamma(C)$. This $C$ is the cost of extracting at $(x_0,y_0)$. This is demonstrated in figure \ref{fig:demonstrateCost}.

The important term from a modeling perspective is the coefficient of $\abs{\nabla\phi}$. The denominator of this term has two summands, $1/v(x,y)$ and $\alpha\psi(x,y)B(x_0,y_0)$. For the first summand, $v(x,y)$ is the walking speed at $(x,y)$ which is determined by the terrain at $(x,y)$. It is calculated by using the steepness at $(x,y)$ and a walking speed function adapted from \citep{IC2017} (discussed in section \ref{sec:walk}). Increasing $v$ makes the level sets move faster, and hence decreases the cost to reach a given point. The second summand, $\alpha\psi(x,y)B(x_0,y_0)$, is a term that penalizes extractors for spending time in highly-patrolled areas. The patrol density function $\psi$ measures the likelihood of an extractor getting caught at position $(x,y)$, and $B(x_0,y_0)$ is the benefit acquired at $ (x_0,y_0) $. If the extractor gets caught, they lose all of the accrued benefit, so the impact of being caught should scale with the benefit gained. When $\psi(x,y)B(x_0,y_0)$ is larger, the level sets move more slowly, and the cost associated with extracting at a point increases.

To compare the summands we include a parameter $\alpha$, which can be thought of as a risk aversion factor. This parameter weighs the willingness of an extractor to assume risk of being captured versus the willingess to expend time or energy. A large $\alpha$ signifies that the extractor will accept less risk of being captured, and will be willing to exert more effort to reach their chosen extraction point. Alternatively, a small $\alpha$ means the extractor is willing to accept a higher risk of being captured to decrease the physical effort required to travel in and out of the protected region. Note that this summand includes information about the extraction point, through the benefit $B(x_0,y_0)$. This means that the cost function must be calculated separately for each extraction point.

\subsubsection{The walking speed equation}\label{sec:walk}

\begin{figure}[h]
\centering
\includegraphics[width=0.75\textwidth]{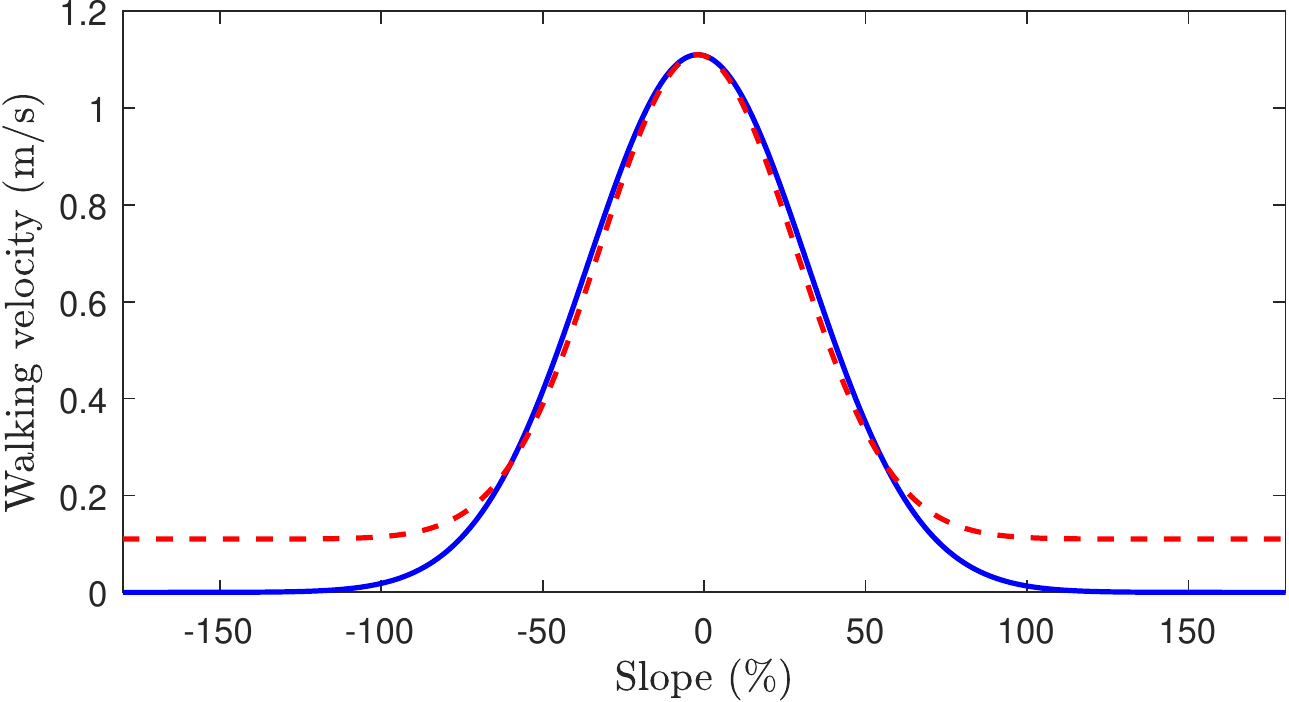}
\caption{Our velocity function \ref{eq:ourIC} (solid blue line), and the \citeauthor{IC2017} function \ref{eq:IC} (dashed red line). Our velocity function decays to zero for high slopes while the function suggest by \citet{IC2017} does not.}
\label{fig:velocityFunctions}
\end{figure}

\Citet{Tobler1993three} gave a formula for walking speed based on slope. Recently, \citet{IC2017} gave an improved formula that they claimed more accurately predicted travel times for humans walking on roads, 
\begin{equation} \label{eq:IC}
v = 0.11+\exp\left(-\frac{\left(100s+2\right)^2}{1800}\right),
\end{equation}
where $100s$ is the grade in percent, and $v$ is the corresponding speed in $\textrm{m}/\textrm{s}$. This formula agrees well with experimental results given in \citep{IC2017}, but has a drawback for our model, that the speed does not vanish as the slope goes to infinity. We hence modify the given equation so that the speed drops to zero as the slope becomes more extreme and use
\begin{equation}\label{eq:ourIC}
v = 1.11\exp\left(-\frac{\left(100s+2\right)^2}{2345}\right),
\end{equation}
which matches the maximum speed, gives similar results for grades less than $60 \%$ and decays to zero for more extreme slopes as can be seen in figure \ref{fig:velocityFunctions}. 

We note that the specific walking velocity function is not crucial in our model, and other functions could be used. We also are able to incorporate areas where the walking speed is faster or slower which could model, for example, overgrown scrub or walking trails. It is also possible to set $v$ to zero in places to model impassable obstacles.

\subsection{The Expected Cost Algorithm}\label{sec:costAlgorithm}
\begin{figure}
	\centering
	\begin{subfigure}[t]{0.31\textwidth}
		\includegraphics[width=\textwidth, trim = 25 10 50 20, clip]{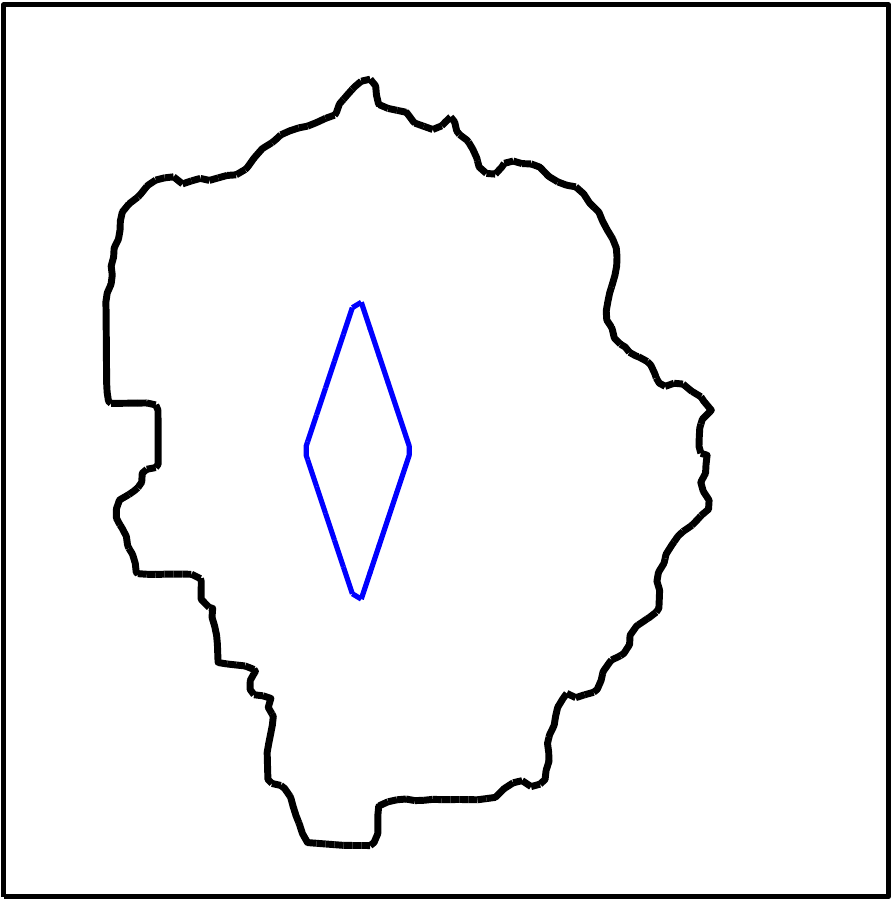} 
    	\caption{Begin by choosing a contour (blue) of fixed benefit.}
    	\label{fig:profit1}
	\end{subfigure}
	~
	\begin{subfigure}[t]{0.31\textwidth}
		\includegraphics[width=\textwidth, trim = 25 10 50 20, clip]{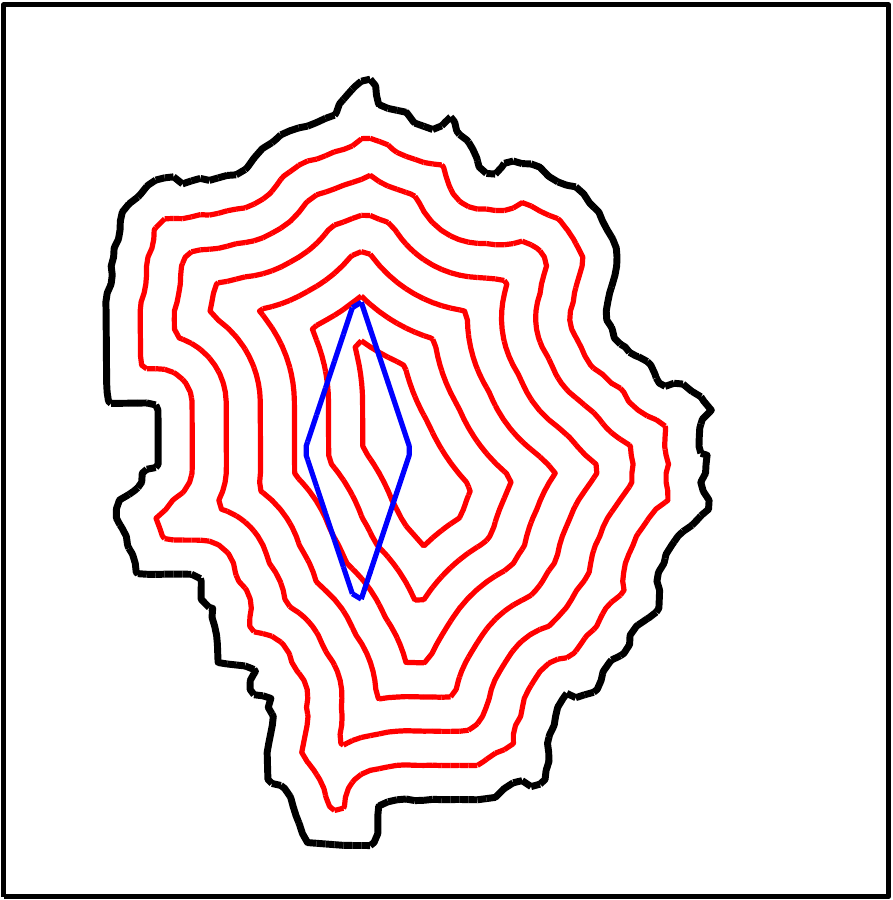}
	    \caption{For this benefit value, evolve cost contours (red) in from the boundary.}
	    \label{fig:profit2}
	\end{subfigure}
	~
	\begin{subfigure}[t]{0.31\textwidth}
		\includegraphics[width=\textwidth, trim = 25 10 50 20, clip]{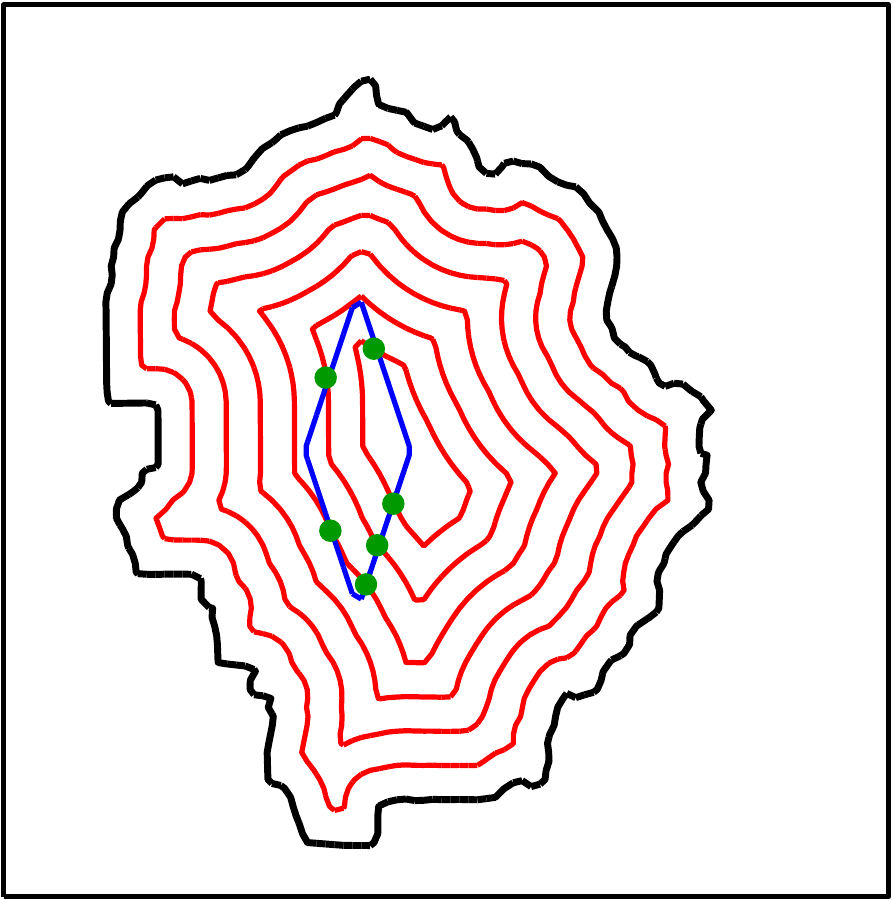}
	    \caption{At intersection points (green), one knows the benefit and cost and thus can calculate the profit.}
	    \label{fig:profit3}
	\end{subfigure}
	\caption{Determining the profit function throughout the region. These images correspond to steps 2, 3 and 4 in the expected cost algorithm (section\ref{sec:costAlgorithm}) respectively.}
\end{figure}

The expected profit function is calculated using the following algorithm.
\begin{enumerate}
\item Since the cost level set equation \ref{eq:costLevelSet} depends on the particular benefit level where extraction takes place, we must solve it for all possible values of $B$. We pick $N$ benefit values in $[\min_\Omega B,\max_\Omega B]$ and evaluate the cost function using each benefit value $B_i$, by performing steps 2-4 for each $B_i$.
\item For a given $B_i$, first find the contour(s) of points where $B(x,y)=B_i$. In figure \ref{fig:profit1}, the black contour is the boundary of the domain, and the blue contour is an equal benefit contour with the same benefit $B_i$. 

\item Evolve the expected cost level set function with the benefit $B_i$ using equation \ref{eq:costLevelSet}. In figure \ref{fig:profit2}, the red contours are the equal-cost level set contours. Each of the contours represents a different cost value.

\item The cost level sets only apply to the points whose benefit has been used to calculate the level sets. That is, we have calculated the cost assuming the benefit gained is $B_i$, so our results are only valid at the points $(x,y)$ such that $B(x,y)=B_i$. Therefore we find the intersections between the cost level contours $\phi(x,y,C)=0$ and the benefit contour $B(x,y)=B_i$. At these intersections (marked in green), we can calculate $C(x,y)$.
\item Repeat steps 2 to 4 for each benefit level $B_i$, until values of $C(x,y)$ are known throughout the region.
\end{enumerate}

Once the expected cost has been found it can be interpolated, and can be used to find the expected profit $P(x,y)=B(x,y)-C(x,y)$. This expression for the profit includes the effect of the defenders.

\subsection{Finding the Pristine Region}
Once the expected profit function is found, the next step is to find the pristine region, the region that extractors will never occupy. In \citet{Albers}, the pristine region was simply the part of the circular forest between its center and the ring of maximum profit. In general, however, the pristine region must be defined differently. We assume that extractors will extract anywhere where the profit is close to the maximum possible profit. Denote $P_{max}$ as the maximum expected profit, and $(1-\varepsilon)P_{max}$ as the lowest expected profit that the extractor will accept for some tolerance $\varepsilon$. We also assume that the extractors will take the optimal path away from their extraction point to the boundary of the region.

The pristine region is determined using the following algorithm.
\begin{enumerate}
\item Find the equal profit contour corresponding to $(1 - \varepsilon)P_{max}$. The extractors will only extract within this high-profit region.
\item Pick points inside the high-profit region uniformly at random, and calculate the optimal path from that point to the boundary of the region. The optimal path can be found easily by solving the gradient descent equation
\begin{equation}
\dd{\vec x}{t}=-\nabla C_{B(\vec x_0)}(\vec x),
\end{equation}
where $\vec x(0)=\vec x_0$ and $C_{B(\vec x_0)}(\vec x)$ is the cost function associated with the benefit value at the extraction point $\vec x_0$.
\item The non-pristine region includes the high-profit region, and all points within some (small) distance of one of the optimal paths found in the previous step. The rest of the region is pristine, and is not traversed by extractors.
\end{enumerate}

\subsection{Metrics for Measuring Patrol Effectiveness}

To measure the effectiveness of a given patrol strategy $\psi$, we use two metrics, although others could be used. A simple measure, as used in \citet{Albers}, is to calculate the proportion of the protected area that is pristine, using,
\begin{equation*}
 \textrm{Pristine Proportion} = \frac{\int_{\Omega} \chi(x,y)\,\mathrm d x\,\mathrm d y}{\int_\Omega\,\mathrm dx\,\mathrm dy},
\end{equation*}
where $\chi(x,y)$ is an indicator function which is 1 in the pristine region and 0 in the non-pristine region.

Some parts of the protected region may be more valuable to protect than other regions, and we could instead choose to measure the proportion of the net value protected using the following integral,
 \begin{equation*}
 \textrm{Proportion of Value Protected} = \frac{\int_{\Omega} B(x,y)\chi(x,y) \,\mathrm d x\,\mathrm d y}{\int_\Omega B(x,y)\,\mathrm dx\,\mathrm dy},
\end{equation*}
where $\chi$ is defined as previously. In place of the benefit function, other functions could be used that measure the relative importance of protecting certain areas.

\subsection{Equivalence of this model and \citep{Albers}}\label{sec:equiv}

Although our model is more sophisticated than that given in \citep{Albers}, it does give the same profit function in some simple symmetric cases. Consider a circular forest of radius 1 with walking speed 1 and a uniform patrol strategy with patrol budget 1 (so $\psi=1/\pi$). Further assume that the benefit associated with extracting at depth $d$ is $B(d)=2d$ and that the cost is $C(d)=d$ (equivalent to the travel time with walking speed 1 in our model). The capture probability associated with extracting at depth $d$ is $\Psi(d)= \int^d_0 \psi(r) dr = d/\pi$, and so the profit function based on equation \ref{eq:AlbersProfit} is 
\begin{equation} \label{eq:AlbersProfitCircle}P(d) = \left(1-\frac{d}{\pi}\right)2d-d=d-\frac{2d^2}{\pi}. \end{equation}

We can arrive at the same result using our model. We assume that the protected region is the unit circle centered at the origin so $\phi_0(x) = 1 - \abs{x}.$ For $x$ in the unit circle, set $B(x) = 2(1-\abs{x})$; this is analogous to setting $B(d) = 2d$ in the Albers model. Further, consider constant walking velocity $v = 1$, constant patrol density $\psi = 1/\pi$ and risk factor $\alpha = 1$. Fix $x_0$ in the unit circle. Note that $\abs{\nabla \phi} = 1$ almost everywhere and the normal velocity $\left(1+2(1-\abs{x_0})/\pi\right)^{-1}$ in equation \ref{eq:costLevelSet} is constant with respect to $x$. Thus the solution to equation  \ref{eq:costLevelSet} is given by \begin{equation} \phi(x,t) = 1 - \abs{x}  - \frac{t}{1 + \frac{2(1-\abs{x_0})}{\pi}}.\end{equation}  Recall that the cost $C(x_0)$ associated with extracting at the point $x_0$ is given implicitly by $\phi(x_0,C(x_0)) = 0$. Solving this equation, we have \begin{equation} C(x_0) = (1-\abs{x_0})\left(1 + \frac{2(1-\abs{x_0})}{\pi} \right) \end{equation} and the profit $P(x_0) = B(x_0) - C(x_0)$ is given by \begin{equation}P(x_0) = (1- \abs{x_0}) \left(1 -  \frac{2(1-\abs{x_0})}{\pi}\right). \end{equation} Note that this profit function is radial and that for any $x_0$ at depth $d$, we have $P(x_0) = d\left(1 - \frac{2d}{\pi}\right)$ exactly as in equation \ref{eq:AlbersProfitCircle}.

In general, with arbitrary benefit function $B(d)$, linear cost $C(d)=cd$, and homogenous patrol density $\psi(d)=p$, Albers gives
\begin{equation} P(d)=(1-pd)B(d)-cd, \end{equation}
while our proposed model gives (with speed set to $v=1/c$ and $\alpha=1$)
\begin{equation} C(d) = d\left(c+pB(d)\right), \end{equation}
and so
\begin{equation} P(d) = B(d)-cd-pdB(d)=(1-pd)B(d)-cd. \end{equation}

The two models agree even if the patrol density is piecewise constant and radially symmetric. If the cost function and/or patrol density are more complicated, then $\phi_t$ will not be constant, and the PDE solution will be more complicated. Our approach to calculating the cost and benefit functions is therefore equivalent to the approach given by \citet{Albers} in simple scenarios but is capable of modeling much more general cases. 

\subsection{Numerical Implementation}

We briefly discuss the numerical implementation of our model. Numerical solutions to Hamilton-Jacobi equations have been well-studied for several years. Recent approaches overcome the curse of dimensionality by using Hopf-Lax formulations and optimization \citep{YTChow}, but since we are concerned with two spatial dimensions, grid-based finite difference methods are sufficient. Recall, a Hamilton-Jacobi equation is an equation of the form $\phi_t + H(x,\nabla \phi) = 0$. Due to the nonlinear dependence on $\nabla \phi$, we cannot use basic differencing methods since they may lead to instabilities or fail to revolve the viscosity solution to the equation \citep{CrandallLions}. Rather, in order to minimize oscillation and track the direction in which information is traveling, we replace the Hamiltonian $H(x,y,\phi_x,\phi_y)$ with a numerical Hamiltonian $\hat H(x,y,\phi_x^+, \phi_x^-, \phi_y^+, \phi_y^-)$ where $\phi_x^+,\phi_x^-$ are the forward and backward difference approximations to $\phi_x$ respectively, and similarly for $\phi_y^+, \phi_y^-$. There are many different acceptable numerical Hamiltonians \citep{OsherShu1991}. We use the Godunov Hamiltonian given by \begin{equation}\label{eq:godunovHamil}\hat H(x,y, \dphi{x}{+},\dphi{x}{-},\dphi{y}{+},\dphi{y}{-}) = \ext{u \in I(\dphi{x}{-},\dphi{x}{+})}\,\,\, \ext{v \in I(\dphi{y}{-},\dphi{y}{+})} H(x,y,u,v) \end{equation} where \begin{equation}I(a,b) = [\min(a,b), \max(a,b)]\end{equation} and \begin{equation} \ext{x \in I(a,b)} = \left \{\begin{matrix} \min_{a \le x \le b} & \text{if } a \le b, \\ \max_{b \le x \le a} & \text{if } a > b. \end{matrix} \right. \end{equation} Performing the minimization or maximization can be difficult if the Hamiltonian is complicated, but in simple cases, they can be resolved explicitly. In our case, the Hamiltonian is $H(x,y,\phi_x,\phi_y) = \tilde v(x,y) \abs{\nabla \phi}$ where the velocity function $\tilde v$ is positive. When the intervals $I(\dphi{x}{-},\dphi{x}{+})$ and $I(\dphi{y}{-},\dphi{y}{+})$ do not include $0$, this Hamiltonian is monotone in each of $\phi_x$ and $\phi_y$ on the rectangle $I(\dphi{x}{-},\dphi{x}{+}) \times I(\dphi{y}{-},\dphi{y}{+})$ and the extrema in \ref{eq:godunovHamil} must occur at the corners of the rectangle. If these intervals do contain $0$, some adjustment is necessary, but one can still explicitly resolve the extrema and find that \begin{equation} \label{eq:GodunovHamiltonianExplicit}\hat H(x,y,\phi^+_x, \phi_x^-,\phi^+_y,\phi_y^-) = \tilde v(x,y) \sqrt{\max\{(\phi_x^-)_+^2,(\phi_x^+)_-^2\} + \max\{(\phi_y^-)_+^2,(\phi_y^+)_-^2\} } \end{equation} where $(A)_+ = \max(A,0)$ and $(A)_- = \min(A,0)$. The Godunov Hamiltonian $\hat H$ gives a first-order approximation to the Hamiltonian $H$. Following \citet{OsherShu1991}, we approximate the derivatives $\phi_x,\phi_y$ at second order and we use second-order total variation diminshing Runge-Kutta time stepping to evolve the solution. In doing so, we have constructed a second order accurate essentially non-oscillatory scheme for \ref{eq:costLevelSet}.

There is an implementation issue which is worth mentioning. Note that the initial function $\phi(x,0)$ can be taken to be precisely the signed distance from $x$ to the initial contour: 
\begin{equation}\phi(x,0) =\text{dist}(x,\Gamma) =  \left\{\begin{matrix} \inf_{y \in \Gamma} \abs{x-y}, & x \text{ inside } \Gamma, \\ -\inf_{y \in\Gamma} \abs{x-y}, & x \text{ outside } \Gamma. \end{matrix} \right. \end{equation}
This is desirable since the distance functions has gradient $1$ almost everywhere and thus if we can preserve this property (that is, if we can ensure that $\phi(x,t)$ is the signed distance to $\Gamma(t)$ for all $t > 0$) we will observe the exact level set motion that we want. However as the level sets evolve, there is some distortion so that numerically for $t > 0$, we no longer have $\phi(x,t) = \text{dist}(x,\Gamma(t))$. If $\abs{\nabla \phi}$ becomes too small, it is difficult to resolve the level sets accurately and if $\abs{\nabla \phi}$ becomes to large, we are required to choose exceedingly fine time discretization or else the numerical solution may develop instabilities. We can prevent these from happening by occasionally replacing $\phi$ with the signed distance function to $\Gamma(t)$. That is, every so often, we halt the time integration, find the current zero level contour $\Gamma(t)$, reset $\phi(x,t) = \text{dist}(x,\Gamma(t))$ and continue. This process is known as \emph{re-distancing} and is discussed by \citet{redistancing}. 

\section{Results}\label{sec:results}
\begin{figure}[htbp]
\centering
\begin{subfigure}[t]{0.49\textwidth}
	{\includegraphics[width=\textwidth,trim=100 40 120 50,clip]{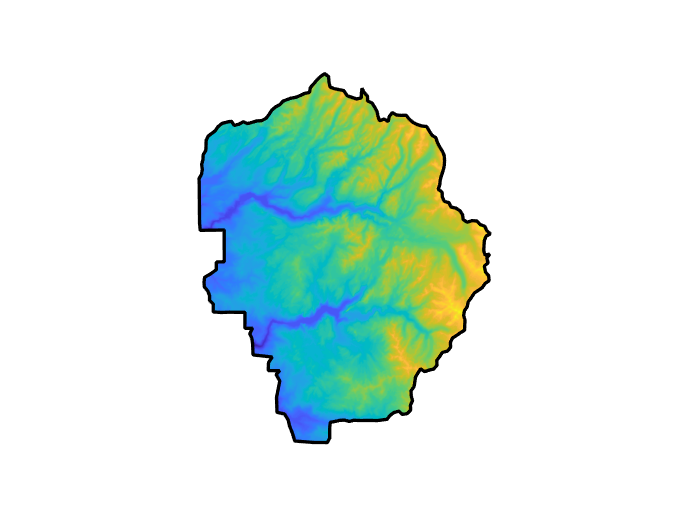}}
	\caption{Yosemite National Park.}
	\label{fig:yosemiteElevation}
\end{subfigure}~
\begin{subfigure}[t]{0.49\textwidth}
	\includegraphics[width=\textwidth,trim=90 80 100 100,clip]{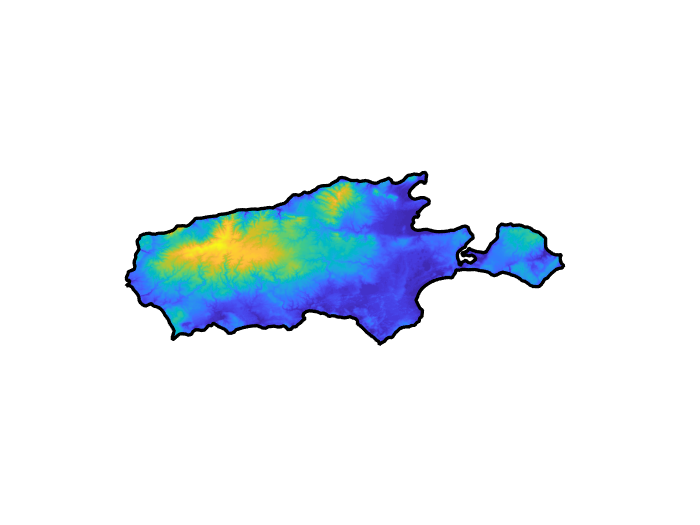}
	\caption{Kangaroo Island.}
	\label{fig:kangarooElevation}
\end{subfigure}

\caption{Elevation profiles for Yosemite National Park and Kangaroo Island. Yellow corresponds to higher elevation, blue to lower elevation (scales different in each figure).}
\label{fig:elevationPlots}
\end{figure}

We present results for our algorithm applied to two real world locations: Yosemite National Park in California, and Kangaroo Island in South Australia. Yosemite National Park is a mountainous area with steep mountains and long valleys. Kangaroo Island has an interesting shape, with a narrow neck separating the main part of the island from a smaller part at the eastern end. For both locations we use real elevation data, sourced from the United States Geographical Survey (Yosemite National Park data) and the Foundation Spatial Data Framework (Kangaroo Island data). The elevation profiles for Yosemite National Park and Kangaroo Island are displayed in figure \ref{fig:elevationPlots}. The data was processed using QGIS \citep{QGIS} and imported to MATLAB using TopoToolbox \citep{topotoolbox2,topotoolbox1}. We apply several patrol strategies identified by \citet{Albers} and \citet{johnson2012patrol} before suggesting some simple and more effective patrols that account for the geometry of the regions. 

\subsection{Yosemite National Park without patrols}
We first consider the case without any patrols, so the cost function depends only on the effort required to travel from the extraction point to the boundary of the park. In figure \ref{:nopatrolling}, we present results for two cases with different benefit functions. Both benefit functions have the same form, a quadratic increase from 0 at the boundary to a maximum value at the point furthest from the boundary, but figure \ref{:nopatrol2} has maximum benefit double that as in figure \ref{:nopatrol}. The lower benefit case has a larger pristine area and smaller high-profit region. In the high-benefit case, there is enough incentive for extractors to expend more effort and extract from more locations within the protected region, obtaining much higher profits (although not doubled). A selection of the optimal paths from extraction points to the boundary of the protected region are also shown. In the sections that follow, we will use $k=8$, the high-benefit case.

\begin{figure}[htbp]
\centering
\begin{subfigure}{0.45\textwidth}
	\includegraphics[width=\textwidth,trim=120 30 150 30,clip]{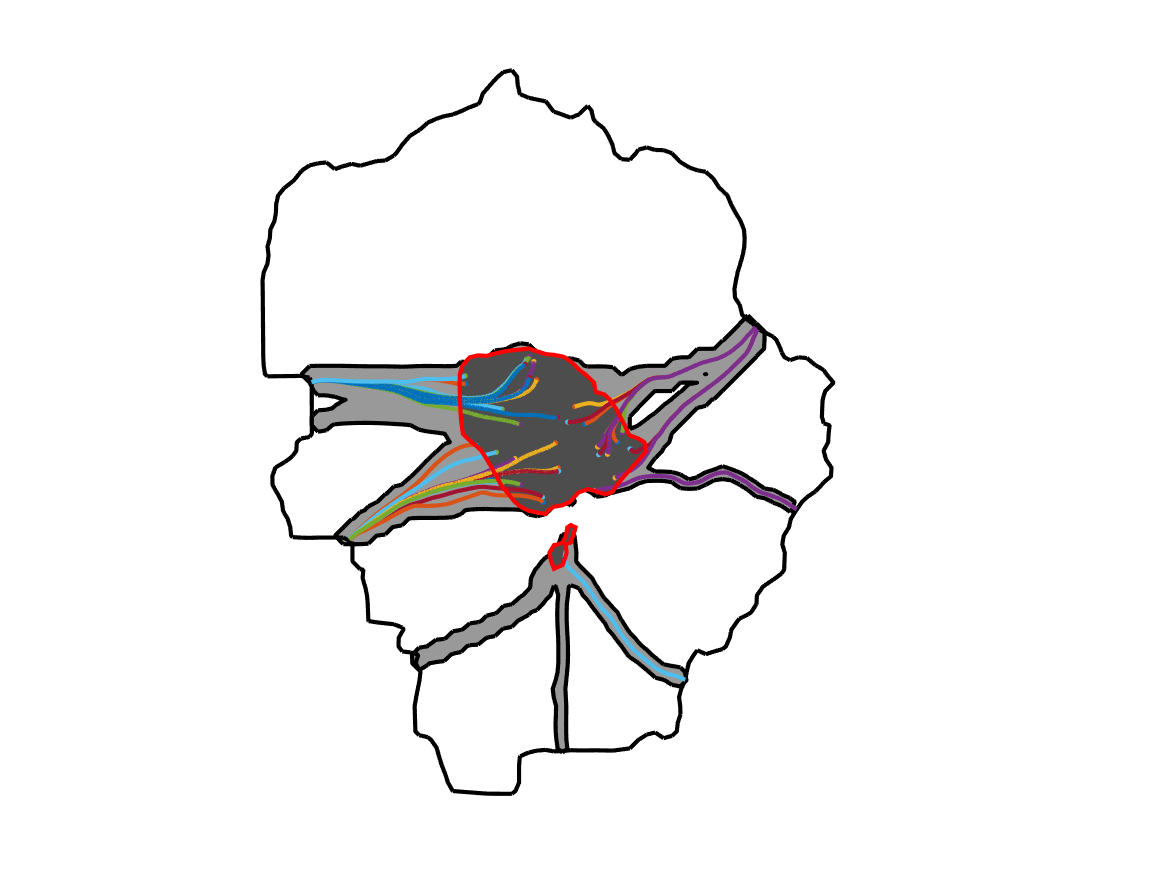}
	\caption{$k=4$, 24 paths out of the total 192 are shown. Maximum profit is $7.06\times10^4$ and pristine proportion is 0.867.}
	\label{fig:nopatrol}
\end{subfigure}~
\begin{subfigure}{0.45\textwidth}
	\includegraphics[width=\textwidth,trim=120 30 150 30,clip]{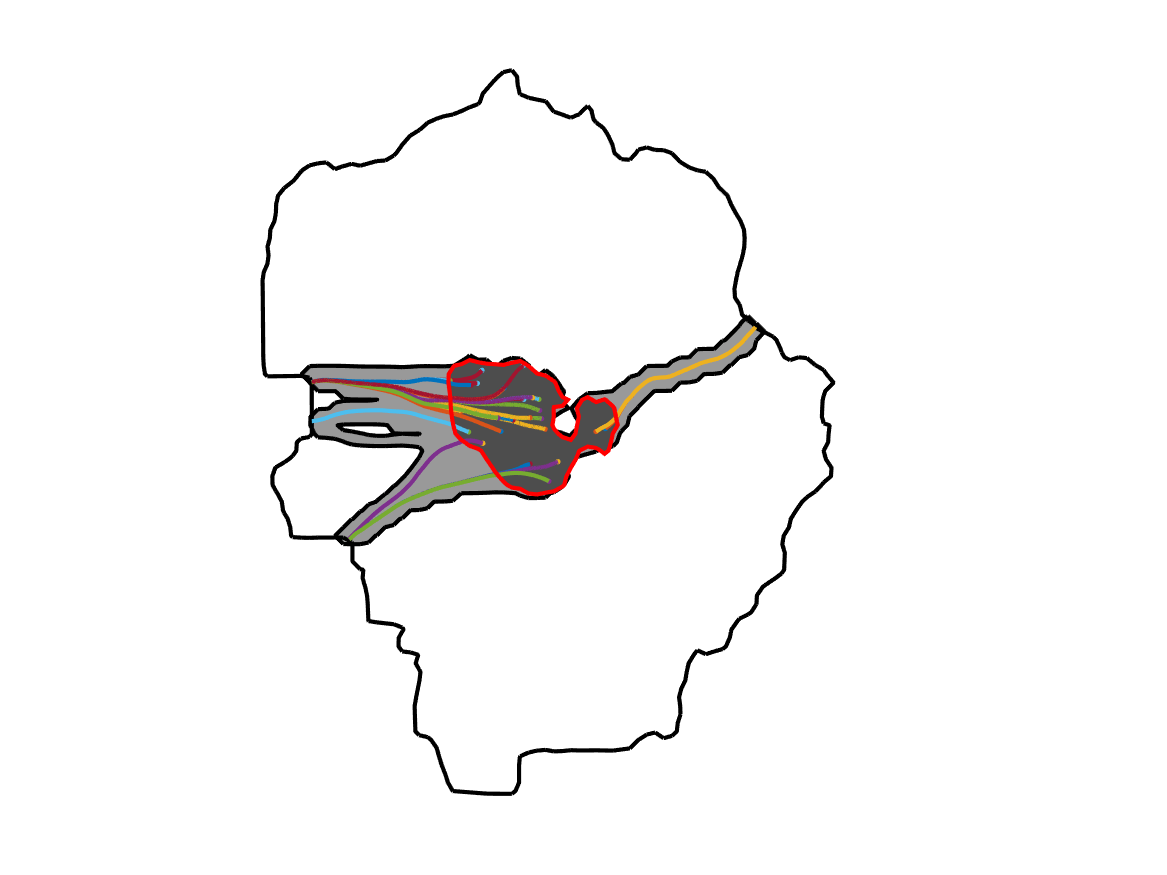}
	\caption{$k=8$, 37 paths out of the total of 294 are shown. Maximum profit is $1.65\times10^5$ and pristine proportion is 0.784.}
	\label{fig:nopatrol2}
\end{subfigure}
\caption{Two cases for Yosemite National Park with no patrols, benefit based on distance from boundary $d$ as $B(d)=kd(2d_m-d)/d_m$ (where $d_m=\max_\Omega d$) and $k=4$ in (a) and 8 in (b). The dark gray denotes the high-profit region where extraction occurs, and light gray shows the envelope of paths the extractors follow when exiting the protected region. Some of the individual paths are shown for illustrative purposes.}
\label{fig:nopatrolling}
\end{figure}

\subsection{Homogeneous patrols}
In figure \ref{:hompatrolling} we consider the simplest nonzero patrol strategy, a homogeneous patrol in which the entire protected area is patrolled with equal intensity. Two cases are shown, both with the high-benefit case $k=8$ from the previous section, but with different patrol budgets $E$. The patrol strategy is simply $\psi=E/A$ where $E$ is the budget and $A$ is the area of the protected region. We specify $\alpha=1$ as the risk-aversion parameter, as this choice gives agreement between our model and that of \citet{Albers} in simple symmetric case (discussed in section \ref{sec:equiv}). As expected, when the patrol budget increases the pristine proportion increases and the maximum profit obtained by the extractor decreases.

\begin{figure}[htbp]
\centering
\begin{subfigure}{0.45\textwidth}
	\includegraphics[width=\textwidth,trim=90 30 110 20,clip]{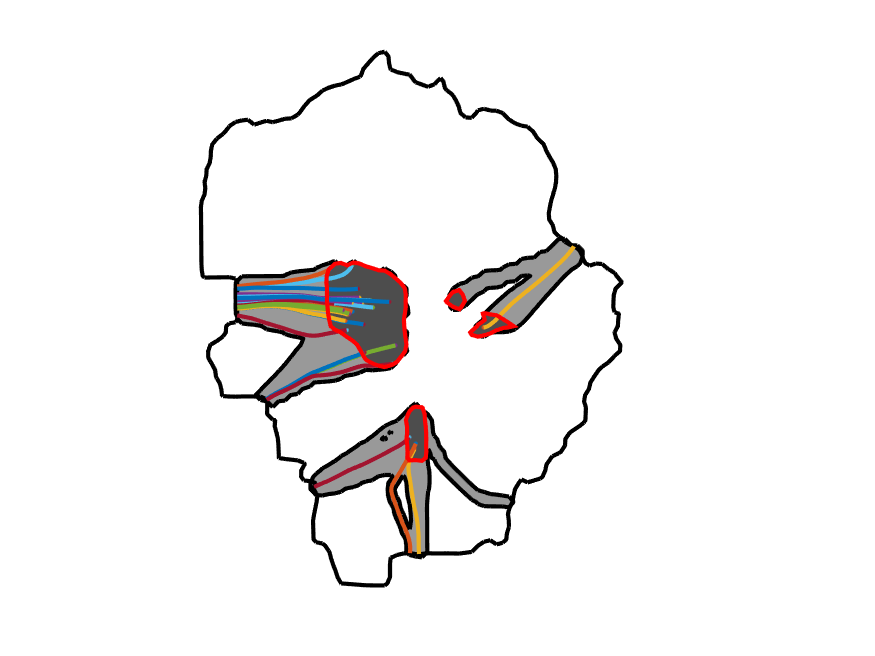} 
	\caption{$E=3\times10^4$, 26 paths out of the total of 205 are shown. Maximum profit is $1.27\times10^5$ and pristine proportion is 0.809.}
	\label{fig:hompatrol}
\end{subfigure}~
\begin{subfigure}{0.45\textwidth}
	\includegraphics[width=\textwidth,trim=90 30 110 20,clip]{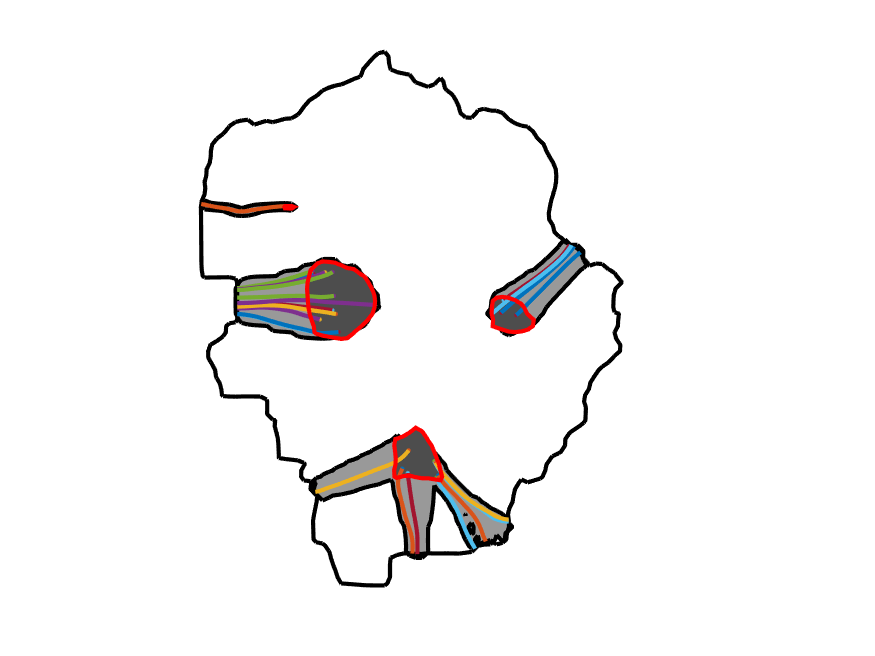} 
	\caption{$E=6\times 10^4$, 23 paths out of the total of 179 are shown. Maximum profit is $9.86\times10^4$ and pristine proportion is 0.841.}
	\label{fig:hompatrol2}
\end{subfigure}
\caption{Two cases for Yosemite National Park with homogeneous patrols, benefit based on distance from boundary $d$ as $B(d)=kd(2d_m-d)/d_m$ (where $d_m=\max_\Omega d$ and $k=8$) and different budgets $E=3\times 10^4$ in (a) and $E=6\times10^4$ in (b).}
\label{fig:hompatrolling}
\end{figure}

In figure \ref{:hompatrolling}, doubling the patrol budget decreases the maximum profit obtained by extractors by roughly 22\%, increases the pristine proportion by roughly 4\%, and increases the proportion of benefit protected by roughly 8.6\%. These are clear improvements, but one perhaps would expect more change upon doubling the patrol density. \Citet{johnson2012patrol} claimed that the homogeneous patrol was not optimal, and this agrees with our intuition that some areas of the protected region will never be penetrated by extractors, so there is no need to patrol there.

\subsection{Band Patrols}
\Citet{johnson2012patrol} identified the band patrol as the optimal strategy in symmetric circular protected regions. In the circular case, the band patrol is a band between two distances from the boundary, $0<d_o<d_i<d_m$, where $d_m$ is the maximum distance from the boundary, with highest patrol density at the outer extent $d_o$, and decreasing density moving towards the inner-most extent of the band $d_i$. In the general case, we set up a band patrol by patrolling between $0.3d_m$ and $0.7d_m$ from the boundary. figure \ref{:bandresults} shows the patrol strategy and the results of our algorithm. We did not implement the algorithm presented by \citet{johnson2012patrol} to find the optimal band patrol, instead testing a number of band patrols and choosing the one that gave the best results.

\begin{figure}[htbp]
\centering
\begin{subfigure}[t]{0.45\textwidth} 
	\includegraphics[width=\textwidth,trim=150 50 160 50,clip]{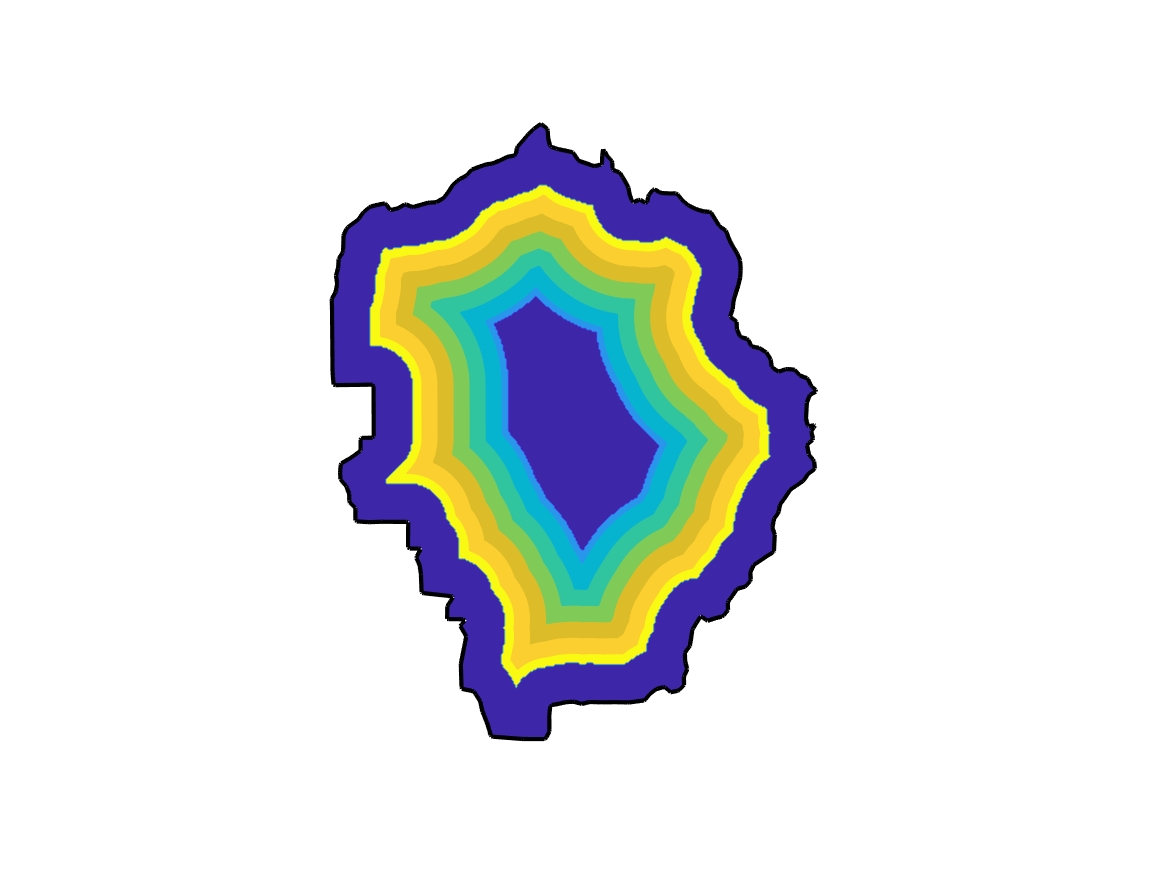} 
	\caption{The patrol strategy. Yellow corresponds to more intense patrolling, blue to no patrolling.}
	\label{fig:bandstrategy}
\end{subfigure}~
\begin{subfigure}[t]{0.45\textwidth}
	\includegraphics[width=.95\textwidth,trim=120 30 160 33,clip]{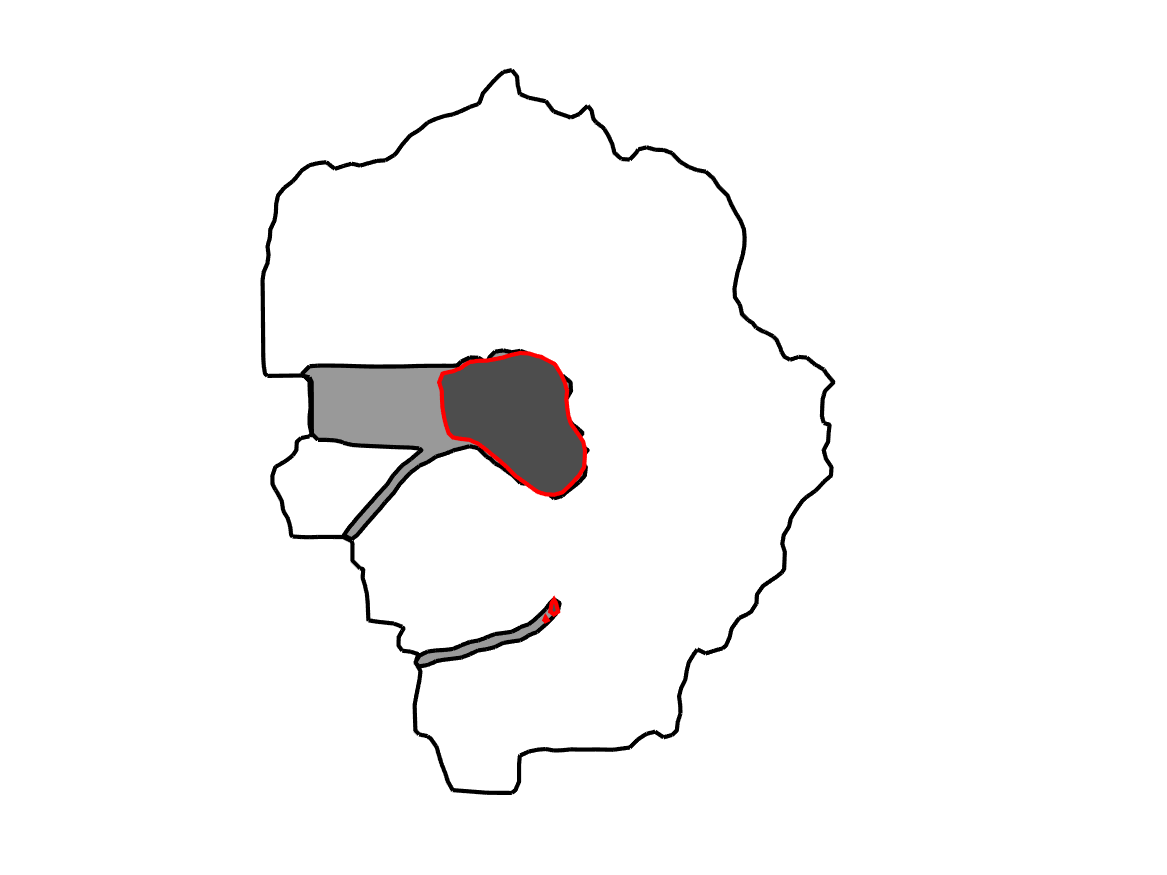}
\caption{The maximum profit is $1.25\times10^5$ and the pristine proportion is 0.901.}
\label{fig:bandpatrol}
\end{subfigure}
\caption{A band patrol in Yosemite National Park. The patrol is based on distance from boundary $d$ and decreases linearly from $d=0.3$ to $d=0.7$, and is zero elsewhere. The benefit is $B(d)=kd(2d_m-d)/d_m$ (where $d_m=\max_\Omega d$ and $k=8$) and the budget $E=3\times 10^4$.}
\label{fig:bandresults}
\end{figure}

The patrol budget in figure \ref{:bandpatrol} is $3\times10^4$, the same as for the homogeneous patrol in figure \ref{:hompatrol}, but the outcome is significantly better for the patrollers. The homogeneous patrol had a pristine proportion 0.809 and was able to protect 76.9\% of the total benefit, whereas the band patrol has a pristine proportion of 0.901 and protects 86\% of the total benefit. The optimal band patrol of \citet{johnson2012patrol} has the property that extractors do not enter the patrolled region, which ours does not. Despite this, it is still superior to the the homogenous patrol (even the homogeneous patrol with twice the budget!) and the other band patrols tested.

\subsection{Asymmetric Patrols}
The patrol strategies employed above were identified by \citet{Albers} and are radially symmetric with respect to the geometry of Yosemite National Park. That is, these strategies depended only on the depth $d$ of the point $(x,y)$: $\psi(x,y) =  \psi^* (d)$. Our method does not require patrols to be radially symmetric, and we here present an asymmetric patrol that outperforms the symmetric patrols discussed previously. Observing the results of the other patrols, it seems that extractors prefer to enter and leave the park at certain portions of the boundary: the portions which are most concave. This is intuitive, as entering at those areas will ensure less travel distance to reach the center. In response, we can design strategies to preferentially patrol those regions through which extractors are more likely to enter. figure \ref{:starresults} shows just such a patrol strategy.

\begin{figure}[htbp]
\centering
\begin{subfigure}[t]{0.45\textwidth} 
	\includegraphics[width=\textwidth,trim=140 50 140 50,clip]{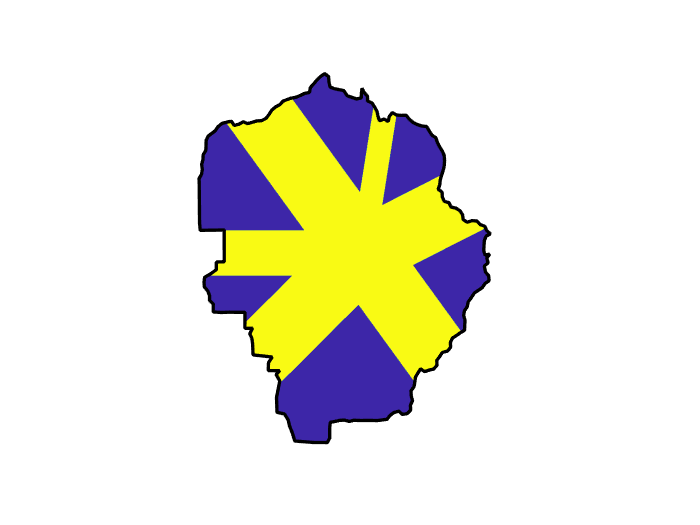} 
	\caption{The patrol strategy. Yellow corresponds to more intense patrolling, blue to no patrolling.}
	\label{fig:starstrategy}
\end{subfigure}~
\begin{subfigure}[t]{0.45\textwidth}
	\includegraphics[width=.95\textwidth,trim=90 30 120 20,clip]{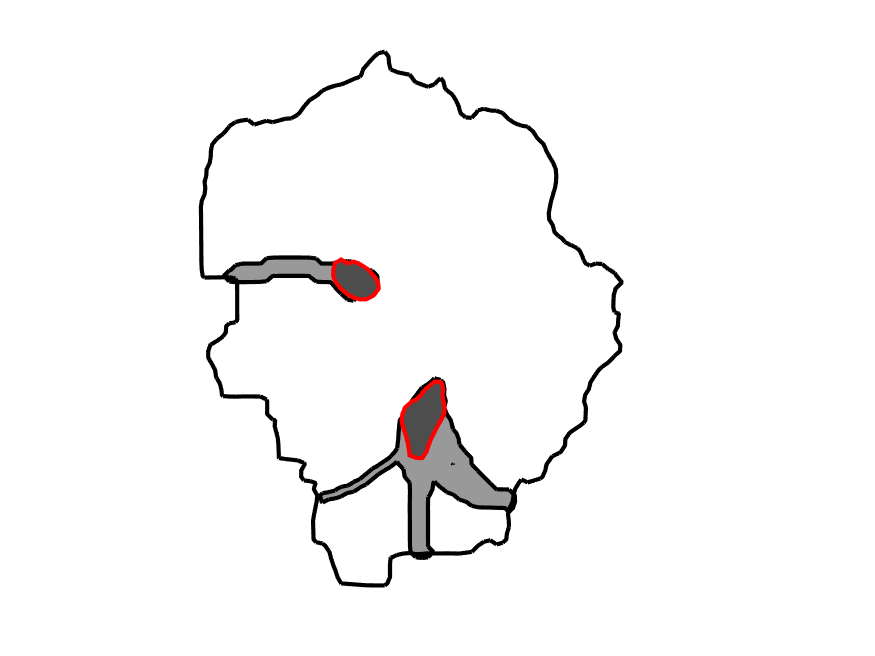}
\caption{The maximum profit is $1.52\times10^5$ and the pristine proportion is 0.920.}
\label{fig:starpatrol}
\end{subfigure}
\caption{A custom patrol in Yosemite National Park which is designed to patrol the concavities in the boundary of the park. The benefit is $B(d)=kd(2d_m-d)/d_m$ (where $d_m=\max_\Omega d$ and $k=8$) and the budget $E=3\times 10^4$.}
\label{fig:starresults}
\end{figure}

Exploiting the geometry of the park, rather than patrolling in a depth-dependent way, the pristine proportion is increased to 0.92 and 90.7\% of the benefit is protected. These results are better than the band or homogeneous patrols with the same budget, and could likely be improved upon even further by finding even more specialized strategies. This result shows the critical importance of the geometry of the protected region. To a very rough approximation Yosemite is fairly circular, but the relatively small-scale concavities are very important in understanding the behaviour of extractors.

\subsection{Kangaroo Island} 
We now apply our method to Kangaroo Island, a small island off the southern coast of Australia. We present two patrol strategies (a homogeneous patrol and a custom designed patrol) which once again emphasize the importance of accounting for the geometry of the protected area. Kangaroo Island is fairly long and narrow with a portion at the eastern end connected to the rest of the island by a small neck.

\begin{figure}[htbp]
\centering
\begin{subfigure}{0.45\textwidth}
	\includegraphics[width=\textwidth,trim=30 85 50 90,clip]{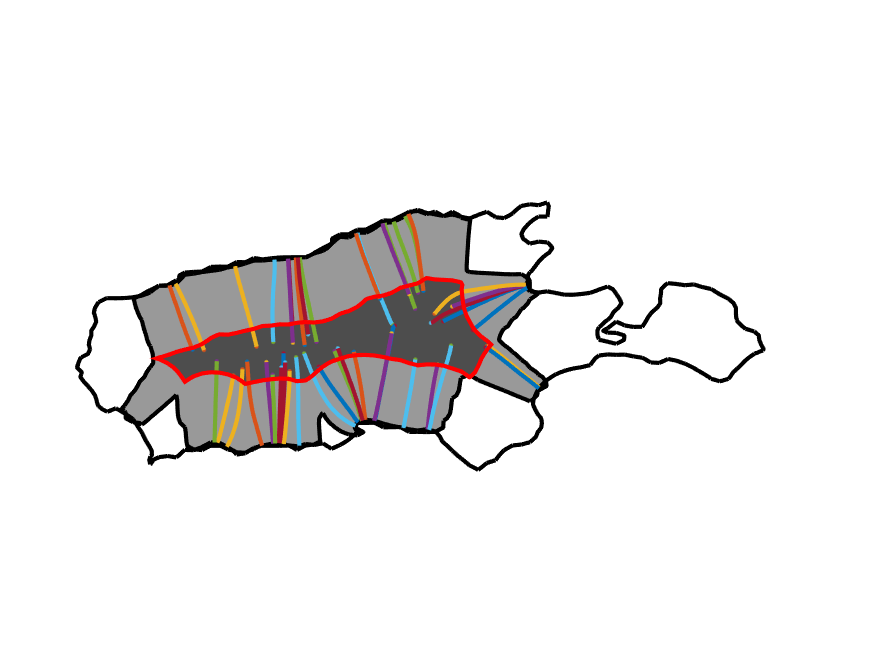} 
	\caption{$E=3\times10^4$, 44 paths out of the total of 431 are shown. Maximum profit is $1.34\times10^5$ and pristine proportion is 0.349.}
	\label{fig:kangahompatrol}
\end{subfigure}~
\begin{subfigure}{0.45\textwidth}
	\includegraphics[width=\textwidth,trim=30 85 50 90,clip]{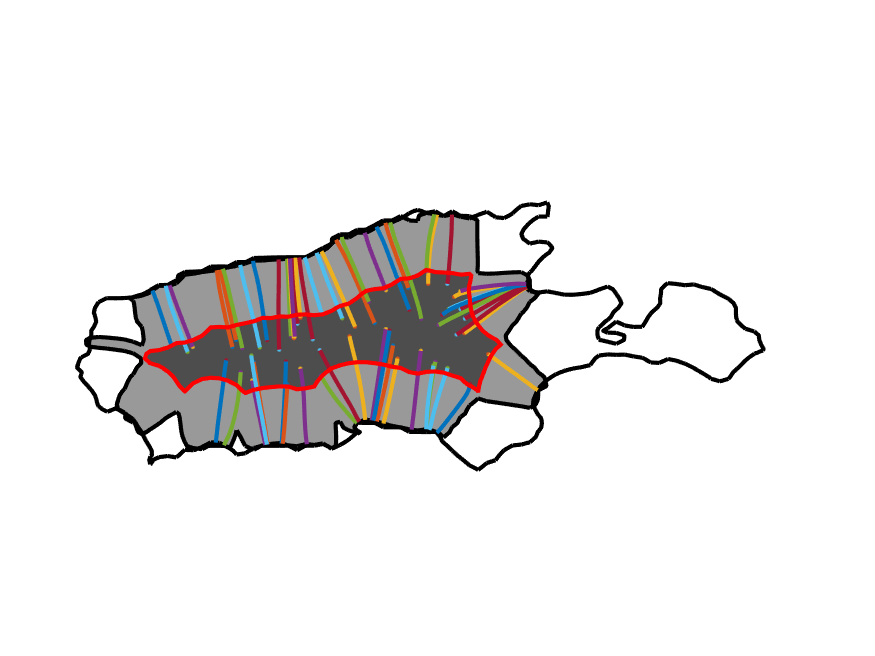} 
	\caption{$E=6\times 10^4$, 58 paths out of the total of 579 are shown. Maximum profit is $1.15\times10^5$ and pristine proportion is 0.305.}
	\label{fig:kangahompatrol2}
\end{subfigure}
\caption{Two cases for Kangaroo Island with homogeneous patrols, benefit based on distance from boundary $d$ as $B(d)=kd(2d_m-d)/d_m$ (where $d_m=\max_\Omega d$ and $k=8$) and different budgets $E=3\times 10^4$ in (a) and $E=6\times10^4$ in (b).}
\label{fig:kangahompatrolling}
\end{figure}

Figure \ref{fig:kangahompatrolling} shows results of homogeneous patrols applied to Kangaroo Island. In this case, we observe that, somewhat counter-intuitively, the increased patrol budget results in a smaller pristine area proportion. The reason for this is that increasing the patrol budget decreases the profit function but causes the high profit area to spread out. That is, with increased patrol the surface profile of the profit function will look more like a plateau which has a smaller maximum than in the low budget case, but has a larger near-maximal area. Although the pristine region is smaller, the maximum profit obtained by the extractors is much lower, decreasing from $1.34\times10^5$ to $1.15\times10^5$, so the patrol does have a significant effect on the extractors. We also observe that a large portion of the non-pristine area comprises the paths leaving the region as opposed to the high-profit region itself. This should inform our decision regarding how to more effectively patrol the region. In any radially symmetric patrol, the high-profit region will also be radially symmetric meaning that it will likely occupy the middle of the island and the paths will cover a very large area as they enter and exit from the north and south side of the island. If we can force the high profit region out of the middle of the island, then the paths will instead travel throughout the east and west portions of the island, hopefully occupying a smaller area. Another key geographic feature of the island is the peninsula at the eastern tip of the island. The peninsula is connected to the island by an isthmus thin enough that, for our purposes, the peninsula can almost be considered an independent region. Extractors will likely be uninterested by the peninsular region since it has much less depth and thus offers much less benefit than the main body of the island. Hence any patrol in that region is likely wasted effort.

With the homogeneous results in mind, figure \ref{:kangaresults} shows the results of a more effective patrol, where the middle of the island is patrolled uniformly and the eastern and western ends are not patrolled. This patrol was able to increase the pristine area proportion to 0.875 compared to the homogeneous patrol which gave 0.349 with the same budget. Again, designing this patrol required some simple observations regarding the geometry of the region, and once again shows the importance of explicitly incorporating geographical information into the model.

\begin{figure}[htbp]
\centering
\begin{subfigure}[t]{0.45\textwidth} 
	\includegraphics[width=\textwidth,trim=90 130 90 120,clip]{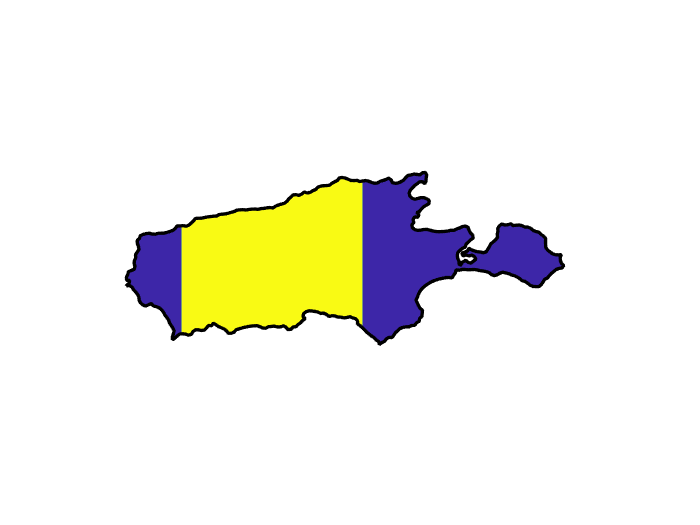} 
	\caption{The patrol strategy. Yellow corresponds to more intense patrolling, blue to no patrolling.}
	\label{fig:kangastrategy}
\end{subfigure}~
\begin{subfigure}[t]{0.45\textwidth}
	\includegraphics[width=.95\textwidth,trim=40 110 60 120,clip]{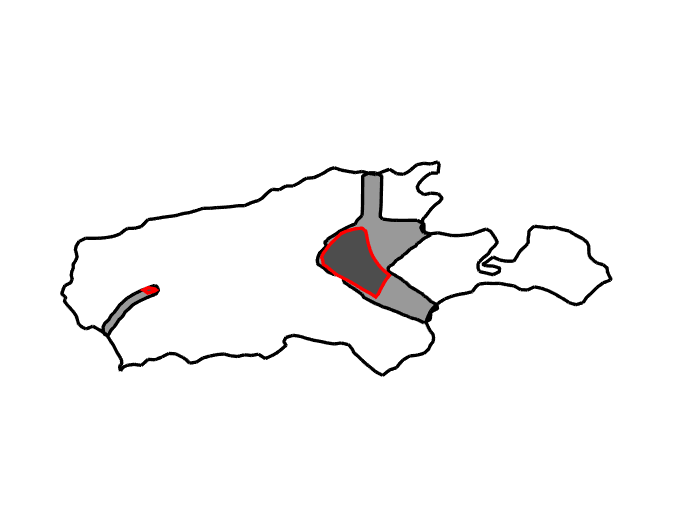}
\caption{The maximum profit is $1.54\times10^5$ and the pristine proportion is 0.875.}
\label{fig:kangapatrol}
\end{subfigure}
\caption{A custom patrol in Kangaroo Island which is designed to push the high profit region away from the center of the island. The benefit is $B(d)=kd(2d_m-d)/d_m$ (where $d_m=\max_\Omega d$ and $k=8$) and the budget $E=3\times 10^4$.}
\label{fig:kangaresults}
\end{figure}

\section{Conclusion}\label{sec:conclusion}

Modeling of environment crime in protected regions has been the subject of previous research, most of which has been focused on discrete network-based methods. Attempts at modelling environmental crime in continuous settings have required restrictive assumptions of symmetry and are only applicable to regions with very simple geometries. In this paper we have formulated a generalized version of the model of \citet{Albers} which can account for regions of arbitrary geometry, does not require any assumptions of symmetry and can incorporate topographical features. In doing so, we have unified environmental crime modeling and continuum, PDE models for tracking human movement.  

The level set method of \citet{osherSethian1988} allows us to perform calculations in regions with complicated geometry. Using the level set method, we track the movement of environmental criminals (extractors) by treating them as a front propagating from the boundary of the protected region toward its center while accounting for changes in travel velocity due to elevation. The level set method is very powerful, and is suitable for many more applications involving movement through regions with complicated geometry. 

Working with arbitrary geometries increases the complexity of the model. Accordingly, we have proposed a novel method for incorporating patrol strategy into the calculation of cost to the extractors and suggested an efficient algorithm for resolving the cost function while performing relatively few level set computations. Further, we have re-defined the pristine area by considering an extraction area (an area of near-maximal profit to extractors) as well as paths that extractors will traverse when leaving the region. This is a necessary deviation from the \citet{Albers} model that arises due to the different topology of the point(s) that maximize the extractors' profit.

We have applied our model to two different regions: Yosemite National Park in California and Kangaroo Island in South Australia. In doing so, we tested several different patrol strategies, including some proposed by \citet{Albers} and \citet{johnson2012patrol}. Then, making some basic observations regarding the geometry of each region, we designed asymmetric patrols which protected the regions more effectively that any of those previously suggested. The success of these asymmetric patrols shows the importance of explicitly treating the geometry. We expect that even more effective patrol strategies exist, and determining them is a potential avenue for further research.

The model presented here could be extended and modified in many ways. For example, the velocity function could be made to depend on ground cover (i.e. higher speeds on open plains and lower speeds in heavy scrub) or other effects such as lakes or paths. Likewise, the benefit and cost functions could be calculated differently depending on the situation. In the proposed model, the cost of travelling to the extraction point is not considered, only the trip out of the region after the crime has been committed. This follows \citet{Albers}, who neglected the cost of the inbound trip. Adding the extra cost of the inbound trip would be a minor modification of the algorithm as presented. We have only considered benefit functions that depend on distance from the boundary, which is likely not realistic. The proposed model can be evaluated with any benefit function (continuity is not necessary), however defining a realistic and accurate benefit function for a given national park would require input from agencies involved in managing the park.

We have defined patrol strategies as density functions, but these do not correspond to patrol routes that would be followed. Determining specific patrol routes that give the desired patrol density is a task requiring significant work. \Citet{fang2017paws} used a discrete model to find explicit patrol routes but in our continuous setting the problem is very complicated. Another opportunity for further research is studying the effects of time-dependent benefit functions. We have used a constant benefit function but in some situations this may not be appropriate. For example, considering animal poaching, one might wish to account for the movement of animals over time. This would require significant modification to the model we have presented.

\section{Acknowledgements} The authors acknowledge the financial support of the NSF (grant DMS-1737770) and the Department of Defense (approved for public release, \#18-666).

\bibliographystyle{plainnat}
\bibliography{bibliography}
\end{document}